\documentclass[10pt,a4paper]{article}%
\usepackage[centertags]{amsmath}
\usepackage{amsfonts}
\usepackage{amssymb}
\usepackage{amsthm}
\usepackage{epsfig}
\usepackage{setspace}
\usepackage{ae}
\usepackage{eucal}
\usepackage[usenames]{color}%
\usepackage{graphicx}

\theoremstyle{plain}
\newtheorem{thm}{Theorem}[section]

\newtheorem{prop}[thm]{Proposition}
\theoremstyle{definition}

\theoremstyle{remark}

\newtheorem{rem}[thm]{Remark}

\begin{document}

\title{ A 3D Euler equation solution with 2D sets of
singularities and data with H\"older continuous
first order derivatives}
\author{J\"org Kampen }
\maketitle

\begin{abstract}
An example of a solution branch of the three dimensional Euler equation Cauchy problem is constructed which develops a singular velocity
component and a singular vorticity component after finite time for some
data which have H\"older continuous first order spatial derivatives. Such a solution branch can be extended beyond a time section at some positive finite
time where a two dimensional set of singularities is located.
\end{abstract}


\section{Introduction}
If $\tau \rightarrow v_i^E(\tau,.), 1\leq i\leq D, \tau\in [0,T)$ denotes a velocity function solution of the $D$-dimensional incompressible Euler equation with final data $v^E_i(T,.)=h_i$, then the function $\tau \rightarrow v_i^{E,-}(\tau,.):=v_i^{E,-}(T-\tau,.), 1\leq i\leq D, t\in [0,T)$ is a solution of the time reversed Euler equation
\begin{equation}
\left\lbrace \begin{array}{ll}
\frac{\partial v^{E,-}_i}{\partial t}-\sum_{j=1}^D v^{E,-}_j\frac{\partial v^{E,-}_i}{\partial x_j}=-\frac{\partial p^{E,-}}{\partial x_i},\\
\\
\sum_{i=1}^D\frac{\partial v^{E,-}_i}{\partial x_i}=0,\\
\\
v^{E,-}_i(0,.)=h_i(.),~\mbox{for}~1\leq i\leq D,
\end{array}\right.
\end{equation}
where $-\frac{\partial p^{E,-}}{\partial x_i}=p_{,i}$
denotes the spatial derivative of the time reversed pressure $p^{E,-}(t,.)=p^E(\tau,.),~\tau=-t$. Here we use classical Einstein notation for
derivatives. The external force term $f_i,~ 1\leq i\leq D$ is assumed to be zero for simplicity, but we shall include it in a short discussion concerning the implications
of our result to Navier Stokes equation with force term in the final section. The
original pressure function $p^E$
is determined by the Euler velocity $v^E_i,~ 1\leq i\leq D$ such that the components of the pressure gradient are 
\begin{equation}
p_{,i}=\sum_{j,m=1}^D\int_{{\mathbb R}^D}K_{D,i}(.-y)\left( \frac{\partial v_m}{\partial x_j}\frac{\partial v_j}{\partial v_m}\right)(t,y)dy,
\end{equation}
where $K_D$ denotes the Laplacian kernel of dimension $D$. For simplicity we consider the special case D = 3, where it is obvious that the following considerations can be applied in the case D > 3 as well. For $\beta_0\in (1.5, 1.5+\alpha_0)$ and $\alpha_0\in \left(0,\frac{1}{2}\right)$ we choose the data 
\begin{equation}
h_1(x)=x_1^{\beta_0}\sin\left(\frac{1}{x_1^{\alpha_0}} \right)\phi_D,~ \mbox{for $x_1\neq 0$, $h_1(0):=0$}
\end{equation}
with weight  function
\begin{equation}
\phi_D=\frac{-2x_D}{\left(1+x_D^2 \right)^3 }\Pi_{i=1}^{D-1}\frac{1}{(1+x_i^2)^2}
\end{equation}
\begin{equation}
\begin{array}{ll}
h_2:=-x_2h_{1,1} = −x_2\left(\beta_0x_1 ^{\beta_0-1}\sin\left(\frac{1}{x_1^{\alpha_0}} \right)  -\alpha_0 x_1^{\beta_0-1-\alpha_0}\cos\left(\frac{1}{x_1^{\alpha_0}} \right)\right) \phi_D\\
\\
-x_2x_1 ^{\beta_0}\sin\left(\frac{1}{x_1^{\alpha_0}} \right)\phi_D\frac{-2x_1}{(1+x_1^2)},~x_1\neq0,~ h_2(0)=0
\end{array}
\end{equation}
and
\begin{equation}
\begin{array}{ll}
\left(\beta_0x_1 ^{\beta_0-1}\sin\left(\frac{1}{x_1^{\alpha_0}} \right) -\alpha_0 x_1^{\beta_0-1-\alpha_0}\cos\left(\frac{1}{x_1^{\alpha_0}} \right)\right)\frac{-x_2^2}{1+x_2^2}\Pi_{i=1}^D\frac{1}{(1+x_i^2)^2}\\
\\
+x_1 ^{\beta_0}\sin\left(\frac{1}{x_1^{\alpha_0}}\right) \frac{4x_1x_2^2}{(1+x_1^2)(1+x_2^2)^2} \Pi_{i=1}^D\frac{1}{(1+x_i^2)^2}~\mbox{for $x_1\neq 0$, and $h_3(0)=0$}.
\end{array}
\end{equation}
Note that $h_i, 1 \leq i \leq 3$ is chosen such that for all $x\in {\mathbb R}^D$ with $x_1\neq $ we have $h_{1,1}(x)+h_{2,2}(x)=x_2h_{1,2}(x)=−h_{3,3}(x)$, or (as $D=3$9
\begin{equation}\label{incompr}
\sum_{i=1}^3h_{i,i}=0~\mbox{ for all $x$ with $x_1\neq 0$}
\end{equation} 
Furthermore, in 7 the singularities of the summands are of the form $\sim x_1^{β_0−(1+α_0)}$ at $x_1=0$ and are related to $h_{1,1}$ and $h_{2,2}$.
 The equation in (\ref{incompr}) holds in $L^2$-sense and in $H^{1,1}$-sense, since the singularities of the second order derivatives are of the form $\sim x^{\beta_0-2(1+\alpha_0)}$
 Indeed as $\beta_0$ becomes close to $1+\alpha_0$
the relation in (\ref{incompr}) holds in
$H^{2-\epsilon}$-sense for small $\epsilon >0$. For our purposes it is relevant that the singularities of $h_{1,1}$
and $h_{2,2}$ are integrable.
We consider a viscosity extension of the equation for the time reversed equation models of the form
\begin{equation}\label{Navleray}
\left\lbrace \begin{array}{ll}
\frac{\partial v^{\nu,-}_i}{\partial t}-\nu\sum_{j=1}^D\frac{\partial^2 v^{\nu,-}_i}{\partial x_j^2} 
-\sum_{j=1}^D v^{\nu,-}_j\frac{\partial v^{\nu,-}_i}{\partial x_j}=\\
\\ \hspace{1cm}-\int_{{\mathbb R}^D}\left( \frac{\partial}{\partial x_i}K_D(x-y)\right) \sum_{j,m=1}^n\left( \frac{\partial v^{\nu,-}_m}{\partial x_j}\frac{\partial v^{\nu,-}_j}{\partial x_m}\right) (t,y)dy,\\
\\
\mathbf{v}^{\nu,-}(0,.)=\mathbf{h},
\end{array}\right.
\end{equation}
to be solved for $\mathbf{v}^{\nu,−}=\left(v^{\nu,-}_1,\cdots,v^{\nu,-}_D\right)^T$ on the domain $[0,\infty)\times {\mathbb R}^D$. Note
that this is not a Navier Stokes equation, because the signs of the nonlinear erms and the sign of the Laplacian are in different relation as in the case of the
Navier Stokes equation. We construct viscosity limits of short time solutions
of the equation in (\ref{Navleray}) with data $h= (h_1,\cdots , h_D)^T$. This viscosity limit turns out to be regular on the time interval $(0,T]$, where $T>0$ is the short time horizon. In the following $H^s$ denotes the standard Sobolev space of order $s\in {\mathbb R}$
and $C^{\delta}$ denotes the space of H\"older continuous functions of H\"older exponent
$\delta\in (0,1)$. Furthermore $C^{1,\delta}$
denotes the space differentiable functions with
H\"older continuous first order derivatives with H\"older exponent $\delta\in (0,1)$. The
main result
\begin{thm}
 Let D = 3. For some time horizon $T>0$ there is a sequence
of regular solutions $v^{\nu_k}_i, 1\leq i\leq  D$ of Cauchy problems with viscosity $(\nu_k)_k$ converging to zero and such that there is a viscosity limit solution 
\begin{equation}
v^{E,−}_i= \lim_{k\uparrow \infty}
v^{\nu_k}_i
\end{equation}
of the time reversed incompressible Euler equation which satisfies $v^{E,−}_i(t,.)\in H^2\cap C^{1,\delta}$
for $t\in (0,T]$   for some $\delta\in (0,1)$. The function $\tau \rightarrow v^E_i(τ, .) := v^{E,−}(t,.)$ with $\tau = T − t \in
[0,T)$ then satisfies the original Euler equation on the time interval $[0,T)$, where
\begin{equation}
v^E_i(T,.)=h,~ 1\leq i\leq D.
\end{equation}
This solution of the Euler equation develops a two dimensional set of singularities after finite time from regular data in $H^2\cap C^{1,\delta}$.
\end{thm}

\begin{rem}
 The construction of singularities can also be done via the time
reversed equation with the data 
\begin{equation}
f_1(x) = x_2x_3
\frac{r^{\beta_0}\sin\left(r^{-\alpha_0}\right)}{(1+r^2)^2},~
f_1(x) =-\frac{1}{2} x_1x_3
\frac{r^{\beta_0}\sin\left(r^{-\alpha_0}\right)}{(1+r^2)^2}
\end{equation}
and
\begin{equation}
f_3(x)=-\frac{1}{2}x_1x_2
\frac{r^{\beta_0}\sin\left(r^{-\alpha_0}\right)}{(1+r^2)^2,~ r=\sqrt{x_1^2+x_2^2+x_3^3}}
\end{equation}
where the choice $1+\alpha_0>\beta_0 > 1.5$, and $\alpha_0\in(0, 0.5)$ is sufficient. In this case
stronger regularity of the data can be obtained for a solution which develops a
singularity at some point in space time after finite time.
\end{rem}

\section{Proof of theorem  }

For given $\nu > 0$ let $G_{\nu}$
denote the fundamental solution of the equation $p_{,t}-\nu \Delta p = 0$. This Gaussian is used for classical representations of solution schemes.
In addition we use the fundamental solution $G_{\epsilon}$ of $p_{,t}-\epsilon \Delta p = 0$ for $\nu > 0$ for
smoothing of the nonlinear terms in the iteration scheme for the first iteration
step. The scheme is designed such that the limit $\epsilon \downarrow 0$ leads to a local time
solution of (\ref{Navleray}) and such that for some subsequence $(\nu_m)_m$ converging to zero
the limit leads to a local time solution of the time reversed Euler equation. For
$\epsilon , \nu >0$ we consider a iterative solution schemes $v^{\nu,\epsilon,-,k}_i,~ 1\leq i\leq D,~ k\geq 0$ related to the fixed point representation in (\ref{Navleray}). We design an iteration scheme
where the convoluted data of the Burgers term and of the Leray projection
term are Lipschitz such that a well-defined viscosity limit can be obtained via
compactness arguments. The choice of the initial data and the design of the first two iteration steps are special. For $k\geq 3$ we have a usual fixed point iteration
scheme related to (\ref{Navleray}). For simplicity we consider the case $D=3$. More precisely
for $k=0$ we define for all $t\geq 0,~1\leq i \leq D$, all $x\in {\mathbb R}^D$
\begin{equation}
v^{\nu,\epsilon,-,(0)}_i(t,x):=\int_{{\mathbb R}^D}
h_i(y)G_{\nu}(t,x;0,y)dy=h_i\ast G_{\nu}
\end{equation}
and for all $1\leq i,j\leq D$ let
\begin{equation}
v^{\nu,\epsilon,-,(0)}_{i,j}(t,x):=\int_{{\mathbb R}^D}
h_i(y)G_{\nu , i}(t,x;0,y)dy=h_i\ast G_{\nu ,j}
\end{equation}
In step $k=1$ of the recursive scheme we may use incompressibility such that we can avoid the data $v^{\nu,\epsilon,-,(0)}_{2,1}$ and
$v^{\nu,\epsilon,-,(0)}_{3,1}$. For $k\geq 1$ we define
\begin{equation}\label{v1}
v^{\nu,\epsilon,-,(k)}=h_i\ast G_{\nu}+B^{(k-1)}_i\ast G_{\nu}+L^{(k-1)}_i\ast G_{\nu}.
\end{equation}

Here we define for $D=3$ the Burgers and Leray projection iteration terms
$B^{(0)}_i, L^{(0)}_i,~1\leq i\leq  D$ such that they are Lipschitz on the time
interval (0, T ] for some time horizon T > 0, and where we define for all $k\geq 0$
\begin{equation}
B^{(k)}_1:=v^{\nu,\epsilon,-,(k)}_1v^{\nu,\epsilon,-,(k)}_{1,1}+v^{\nu,\epsilon,-,(k)}_2v^{\nu,\epsilon,-,(k)}_{1,2}
+  v^{\nu,\epsilon,-,(k)}_3v^{\nu,\epsilon,-,(k)}_{1,3},
\end{equation}
and
\begin{equation}\label{Bk2}
\begin{array}{ll}
B^{(k)}_2:=\left( v^{\nu,\epsilon,-,(k)}_1v^{\nu,\epsilon,-,(k)}_{2}\right)\ast_{sp} G_{\epsilon,1} -\left( v^{\nu,\epsilon,-,(k)}_{1,1}v^{\nu,\epsilon,-,(k)}_{2}\right)\ast_{sp}G_{\epsilon}\\
\\
\hspace{1cm}+v^{\nu,\epsilon,-,(k)}_2v^{\nu,\epsilon,-,(k)}_{2,2}
+ v^{\nu,\epsilon,-,(k)}_3v^{\nu,\epsilon,-,(k)}_{2,3},
\end{array}
\end{equation}
and
\begin{equation}\label{Bk3}
\begin{array}{ll}
B^{(k)}_3:=\left( v^{\nu,\epsilon,-,(k)}_1v^{\nu,\epsilon,-,(k)}_{3}\right)\ast_{sp} G_{\epsilon,1} -\left( v^{\nu,\epsilon,-,(k)}_{1,1}v^{\nu,\epsilon,-,(k)}_{3}\right)\ast_{sp}G_{\epsilon}\\
\\
\hspace{1cm}+v^{\nu,\epsilon,-,(k)}_2v^{\nu,\epsilon,-,(k)}_{3,2}
+ v^{\nu,\epsilon,-,(k)}_3v^{\nu,\epsilon,-,(k)}_{3,3}.
\end{array}
\end{equation}
For the first iteration Leray projection terms we define for all $1 \leq i\leq D$
\begin{equation}\label{Bk4}
\begin{array}{ll}
L^{(k)}_i:= K_{,i}\ast_{sp}
2{\Bigg (}v^{\nu,\epsilon,-,(k)}_{1,1}v^{\nu,\epsilon,-,(k)}_{1,1}+v^{\nu,\epsilon,-,(k)}_{2,2}v^{\nu,\epsilon,-,(k)}_{2,2}\\
\\
+v^{\nu,\epsilon,-,(k)}_{3,3}v^{\nu,\epsilon,-,(k)}_{3,3}
+v^{\nu,\epsilon,-,(k)}_{2,3}v^{\nu,\epsilon,-,(k)}_{3,2}+
\left( v^{\nu,\epsilon,-,(k)}_{1,2}v^{\nu,\epsilon,-,(k)}_{2}\right)\ast_{sp} G_{\epsilon,1}\\
\\
 -\left( v^{\nu,\epsilon,-,(k)}_{1,2,1}v^{\nu,\epsilon,-,(k)}_{2}\right)\ast_{sp} G_{\epsilon}
+\left( v^{\nu,\epsilon,-,(k)}_{1,3}v^{\nu,\epsilon,-,(k)}_{3}\right)\ast_{sp} G_{\epsilon,1}\\
\\
 -\left( v^{\nu,\epsilon,-,(k)}_{1,3,1}v^{\nu,\epsilon,-,(k)}_{2}\right)\ast_{sp}G_{\epsilon}{\Bigg )}.
\end{array}
\end{equation}
Finally, for all $1 \leq i \leq D$ and $k\geq  1$ we define
\begin{equation}
v^{\nu,\epsilon,-,(k)}_{i,j}=h_i\ast_{sp}G_{\nu}+B_i^{(k-1)}\ast G_{\nu,j}+L^{(k-1)}_i\ast G_{\nu,j}
\end{equation}
The iteration scheme above is designed such that such that for the data $h_i, 1\leq i \leq  D$ chosen above we have a well-defined viscosity limit in the sense that there
exists a sequence $(\nu_m)_m$ with $\lim_{m\uparrow \infty}\nu_m = 0$ 
such that $v^{E,−}_i= \lim_{m\uparrow \infty}v^{\nu,−,(k)}_i$ 
is a well defined regular function which satisfies the time-reversed Euler equation on on the time interval $(0,T]$ for some time horizon $T>0$. We are mainly
interested in a construction which implies that we have nondegenerate increments
$\delta v^{\nu,\epsilon,-,(k)}_i:= v^{\nu,\epsilon,-,(k)}_i-h_i\ast_{sp} G_{\nu}$
for $1 \leq i \leq  D$ and all $\nu > 0$ and $k \geq  1$ and nondegenerate derivatives of these
increments. The first observation behind this scheme is that a symmetry of the
first order derivatives of the Gaussian can be combined with Lipschitz continuity
of convoluted Burgers term data $B^{(k)}_i$
and Leray projection term data $L^{(k)}_i$, such
that compactness arguments can be applied to a scheme of the form (21) in order
to obtain a viscosity limit. Symmetry of the first order spatial derivative of the
Gaussian means that for all $y\in {\mathbb R}^D$ we have
\begin{equation}
\left(\frac{-2y_i}{4\nu t} \right)G_{\nu}(t,y)=\left(\frac{-2y^-_i}{4\nu t} \right)G_{\nu}(t,y^-),
\end{equation}
where $y^-=\left(y^-_1,\cdots,y^-_D\right)$ with $y^-_i=-y_i$ and $y_j=y^-_j$ for all $j\in\left\lbrace 1,\cdots,D \right\rbrace \setminus \left\lbrace i\right\rbrace$.
This implies that for any spatially global Lipschitz continuous function $F$ (which
will be Burgers and Leray projection iteration terms $B^{(k-1)}_i$
or $L^{(k-1)}_i$ for our purposes) we have
\begin{equation}
\begin{array}{ll}
{\Big |}F\ast G_{\nu,i}{\Big |}={\Big |}\int F(t-\sigma,x-y)\left(\frac{-2y_i}{4\nu \sigma}\right)  G_{\nu}(\sigma,y)dyd\sigma{\Big |}\\
\\
={\Big |}\int_{y,y_i\geq 0} \left( F(t-\sigma,x-y)-F(t-\sigma,x-y^{-})\right) \left(\frac{-2y_i}{4\nu \sigma}\right)  G_{\nu}(\sigma,y)dyd\sigma{\Big |}\\
\\
\leq L{\Big |}\int_{y,y_i\geq 0} \left(\frac{4y^2_i}{4\nu \sigma}\right)  G_{\nu}(\sigma,y)dyd\sigma{\Big |}=: 4LM_2
\end{array}
\end{equation}
where $L$ is a spatial Lipschitz constant of the function $F$ which is uniform with
respect to time, and $M_2$
is a finite second moment constant. We mention here
that this second moment constant is small in the sense that the growth of the nonlinear Euler terms can be offset by potential damping. However in this paper we concentrate on local time solutions. The second observation is that the first order spatial derivatives of the spatial convolution of the data term $h_i\ast_{sp}G_{\nu,j}$
degenerate as the viscosity tends to zero at time $t>0$. Indeed, for $t>0$ we have for any bounded continuous function $f$ with $\sup_{x\in {\mathbb R}^D} |f (x)| \leq C < \infty$ degeneracy in the viscosity limit as with $y_i= \nu tz_i, 1 \leq i \leq D$ we have
\begin{equation}
\begin{array}{ll}
{\Big |}\int f(x-y)\left(\frac{-2y_i}{4\nu t}\right) \frac{1}{\sqrt{4\nu t}^D} \exp\left(-\frac{|y|^2}{4\nu t} \right) dy\sigma{\Big |}\\
\\
\leq {\Big |}\int Cz_i \frac{1}{\sqrt{4\nu t}^D} \exp\left(-|z|^2\nu t \right)(\nu t)^D dz\sigma{\Big |}\downarrow 0~\mbox{for $\nu\downarrow 0$}
\end{array}
\end{equation}
A third observation that we have polynomial time decay of the increments $\delta v^{\nu,\epsilon,-,(k)}_i(t,.)$ due to the multiplicative structure of the nonlinear terms (convolution with the Gaussian decreases the order of polynomial decay at spatial
infinity a bit but the multiplicative effect is stronger for the data chosen). More
precisely, for the data functions $h_i, 1 \leq i \leq D$ we have for all $m \geq 0$
\begin{equation}
h_i\in {\cal C}^4_{pol,m},~ 1\leq i\leq D,
\end{equation}
where
\begin{equation}
C^l_{pol,m}:=\lbrace f:{\mathbb R}^D\rightarrow{\mathbb R}:\exists c>0~ \forall |x|\geq 1~ \forall 0\leq |\gamma|\leq m~ {\Big |}D^{\gamma}_xf(x){\Big |}\leq \frac{c}{1+|x|^l}\rbrace
\end{equation}
From the multiplicative structure of the nonlinear terms it follows that for $m \geq  0$
\begin{equation}
\delta v^{\nu,\epsilon,−,(k)}_i
(t, .)\in {\cal C}^4_{pol,m}, v^{\nu,\epsilon,−,(k)}_i(t, .)\in C^{4-\mu}_{pol,m}
\end{equation}
for some $\mu \in (0,1)$ which accounts for the convolution effects of the linear
term of the the convoluted initial data. Hence the regularity of the solution
approximations $v^{\nu,-,(k)}_i$ is a matter of local analysis. In order to get a viscosity
limit we need estimates which are independent of $\nu > 0$. Next we consider the
regularity of the iterative approximations $v^{\nu,\epsilon,−,k}_i, 1 ≤ i ≤ D, k ≥ 0$ of the
local time solution functions $v^{\nu,−}_i, 1 \leq i \leq D$ and $v^{E,−}_i, 1 \leq  i \leq D$, which are
obtained by limits of zero subsequences $(\epsilon_l)_l$ and $(\nu_m)_m$ where the iteration index $k\uparrow \infty$. According to
the previous remarks it is essential to understand the regularity of the functional
increments $v^{\nu,\epsilon,-,k}_i, 1 \leq i \leq D, k \geq 0$, especially for the first iteration steps
$k = 1, 2$. The main steps are
\begin{itemize}
\item[i)] For $t>0$ and all $1 \leq i,j \leq  D$ there is a finite constant $C > 0$ which is independent of $\nu$ and $\epsilon >0$ such that for $t> 0$ and all $x\in {\mathbb R}^D$
\begin{equation}
{\big |}\delta v^{\nu,\epsilon,-,(1)}_{i,j}(t,x){\big |}\leq \left( C+C|x_1|^{2\beta_0-3}\right)\Pi_{i=1}^D\frac{1}{(1+x_i^2)^2}
\end{equation}
\item[ii)] For $t>0$ and all $1\leq i \leq D$ and $k\geq 1$ there is a finite constant $C > 0$
which is independent of $\nu$ and $\epsilon>0$ and time such that for all $1 \leq i \leq D$ we have for all $1\leq i\leq D$
\begin{equation}
{\big |}\delta v^{\nu,\epsilon,-,(k)}_i(t,.){\big |}_{H^2\cap C^{1,\delta}}\leq C
\end{equation}
for some $\delta > 0$. Furthermore, there exists $T > 0$ such that for all $t\in [0, T ]$
we have for all $k \geq 1$
\begin{equation}
\max_{1\leq i\leq D}{\big |}\delta v^{\nu,\epsilon,-,(k+1)}_i(t,.){\big |}_{H^2\cap C^{1,\delta}}\leq \frac{1}{2}\max_{1\leq i\leq D}{\big |}\delta v^{\nu,\epsilon,-,(k)}_i(t,.){\big |}_{H^2\cap C^{1,\delta}}.
\end{equation}
\item[iii)] there exists zero sequences $(\epsilon_m)_m$ and $(\nu_m)_m$such that 
\begin{equation}
v^{E,-}_i=\lim_{\nu_m\downarrow 0,\epsilon_m\downarrow 0} h_i\ast_{sp}G_{\nu_m}+\lim_{\nu_m\downarrow 0,\epsilon_m\downarrow 0}\delta v^{\nu_m,\epsilon_m,-,(k+1)},~ 1\leq i\leq D,
\end{equation}
solves the time reversed Euler equation on the time interval $[0,T]$ in $H^1$-sense and classically on the time interval $(0,T]$ and such that for all $t\in (0, T ]$ there exists a finite $C>0$ such that
\begin{equation}
{\big |}v^{E,-}_i(t,.){\big |}_{H^2\cap C^{1,\delta}}\leq C
\end{equation}
for some $\delta >0$.
\end{itemize}

Ad i) we first observe that for $t\geq 0$, all $\nu,\epsilon>0$ and all $x\in {\mathbb R}^D$
\begin{equation}
\begin{array}{ll}
{\big |}v^{\nu,\epsilon,-,(0)}_{1,1}(t,x){\big |}={\big |}h_1\ast_{sp}G_{\nu,1}(t,x){\big |}
={\Bigg |}\left( x_1^{\beta_0}\sin\left(\frac{1}{x_1^{\alpha_0}} \right) \phi_D\ast_{sp} G_{\nu,1}\right) 
(t,x){\Bigg |}\\
\\
\leq 2{\Bigg |}\int_{y,y_1>0}\left( (x_1-y_1)^{\beta_0}\phi_D(x-y) {\Big |}G_{\nu,1}(t,y){\Big |}\right)dy {\Bigg |}\leq Cr^{\beta_0-1}\phi_D
\end{array}
\end{equation}
for some finite constant $C>0$ which is independent of $\nu,\epsilon$ and $t>0$. Analogous upper bounds hold for 
${\Big |}v^{\nu,\epsilon,-,(0)}_j(t,x){\Big |}$ and for ${\Big |}v^{\nu,\epsilon,-,(0)}_{j,j}(t,x){\Big |}$ for $j=2,3$. Furthermore for some finite $C>0$ we have
\begin{equation}
{\big |}v^{\nu,\epsilon,-,(0)}_{j,k}(t,x){\big |}\leq Cr^{\beta_0}\phi_D,~ 2\leq k\leq 3, 1\leq j\leq 3.
\end{equation}
An analogous upper bound holds for ${\big |}v^{\nu,\epsilon,-,(0)}_1(t,x){\big |}$. In step $k=0$ we have no $\nu$-independent upper bound for
\begin{equation}\label{ub}
{\big |}v^{\nu,\epsilon,-,(0)}_{j,1}(t,x){\big |},~ 2\leq j\leq 3.
\end{equation}
Using the product rule for derivatives of multiplication we can avoid the terms
in (\ref{ub}) in a recursion. The first order spatial derivatives of (15) have the representation
\begin{equation}
v^{\nu,\epsilon,-,(1)}_{i,j}=:h_i\ast_{sp} G_{\nu,j}+B^{(0)}_i\ast G_{\nu,j}+L^{(0)}_i\ast G_{\nu,j},
\end{equation}
where $B^{(0)}_i,L^{(0)}_i,~ 1\leq i\leq 3$ are defined above. Using these upper bounds of the data $v^{\nu,\epsilon,-,(0)}_j$ and of the first order spatial
derivatives which enter in the definition of the Burgers and Leray projection
terms of the first approximation (16), (17), (18), and (19) we have global Lipschitz continuity of the Burgers terms and the Leray projection terms, i.e., for all $x,y\in {\mathbb R}^D$ we have for given $t>0$
\begin{equation}
{\big |}B^{(0)}_i(t,x)-B^{(0)}_i(t,y){\big |}\leq L|x-y|,~ 1\leq i\leq D,
\end{equation}
and
\begin{equation}
{\big |}L^{(0)}_i(t,x)-L^{(0)}_i(t,y){\big |}\leq L|x-y|,~ 1\leq i\leq D,
\end{equation}
for some Lipschitz constant L which is independent of ν, ǫ and t > 0. Here we
remark that in the definition of these Burgers and Leray projection terms the
derivatives in 36 are not involved such that we can estimate terms of the form
\begin{equation}
\left(v^{\nu,\epsilon,-,(0)}_1v^{\nu,\epsilon,-,(0)}_2\right)\ast_{sp} G_{\epsilon ,1} 
\end{equation}
in the Burgers term and terms of the form
\begin{equation}
K_{,i}\ast_{sp} 2\left( \left(v_{1,2}^{\nu,\epsilon,-,(0)}v_2^{\nu,\epsilon.-,(0)}\right)\ast G_{\epsilon,1}-v_{1,2,1}^{\nu,\epsilon,-,(0)}v_2^{\nu,\epsilon,-,(0)} \right) 
\end{equation}
in the Leray projection term. In the latter case we note that we have upper bounds  ${\big |}v^{\nu,\epsilon,-,(0)}_{1,2,1}(t,.){\big |}\leq Cr^{\beta_0-1}$ for some  $C$ independent of $\nu,\epsilon$ and $t$. The statement of item i) follows, where we have
\begin{equation}
\delta v^{\nu,\epsilon,-,(1)}_{i,j}(t,.)\in C^{\delta}~ \mbox{for $\delta \in (0,2\beta_0-3)$}.
\end{equation}
Ad ii), the local time contraction result follows from a representation of the
functional increments and its first order spatial derivatives. It is essential to
consider the first order derivatives. For all $1\leq  i \leq  D$ and $k\geq 1$ we have
\begin{equation}\label{deltavk1}
\delta v^{\nu,\epsilon,-,(k)}_{i,j}=\delta B_i^{(k-1)}\ast G_{\nu,j}+\delta L^{(k-1)}_i\ast G_{\nu,j}
\end{equation}
where
\begin{equation}
\begin{array}{ll}
\delta B^{(k-1)}_i:=B^{(k)}_i-B^{(k-1)}_i,~ 
\delta L^{(k-1)}_i:=L^{(k)}_i-L^{(k-1)}_i
\end{array}
\end{equation}
In order to observe how functional  increments can be extracted from interpretations such as in (\ref{deltavk1}) we consider the essential cases based on (\ref{Bk2}) and (\ref{Bk4}) above. Since we consider the contraction strating with $k+1= 2$ we may rewrite the Burgers term for this purpose in a more original form
\begin{equation}\label{Bk2*}
\begin{array}{ll}
B^{(k)}_2:=\left( v^{\nu,\epsilon,-,(k)}_1v^{\nu,\epsilon,-,(k)}_{2,1}\right)\ast_{sp} G_{\epsilon}+v^{\nu,\epsilon,-,(k)}_2v^{\nu,\epsilon,-,(k)}_{2,2}\\
\\
\hspace{1cm}+ v^{\nu,\epsilon,-,(k)}_3v^{\nu,\epsilon,-,(k)}_{2,3}.
\end{array}
\end{equation}
For the some reason we rewrite the Leray projection term, where for all $1 \leq i\leq D$
\begin{equation}\label{Bk4*}
\begin{array}{ll}
L^{(k)}_i:= K_{,i}\ast_{sp}
2{\Bigg (}v^{\nu,\epsilon,-,(k)}_{1,1}v^{\nu,\epsilon,-,(k)}_{1,1}+v^{\nu,\epsilon,-,(k)}_{2,2}v^{\nu,\epsilon,-,(k)}_{2,2}\\
\\
+v^{\nu,\epsilon,-,(k)}_{3,3}v^{\nu,\epsilon,-,(k)}_{3,3}
+v^{\nu,\epsilon,-,(k)}_{2,3}v^{\nu,\epsilon,-,(k)}_{3,2}+
\left( v^{\nu,\epsilon,-,(k)}_{1,2}v^{\nu,\epsilon,-,(k)}_{2,1}\right)\ast_{sp} G_{\epsilon}\\
\\
+\left( v^{\nu,\epsilon,-,(k)}_{1,3}v^{\nu,\epsilon,-,(k)}_{3,1}\right)\ast_{sp} G_{\epsilon}{\Bigg )}.
\end{array}
\end{equation}
Using the uniform upper bounds for $v^{\nu,\epsilon,-,(k-1)}_{i,j}$
for the involved indices of the
previous item, regularity for the Poisson equation, the Young inequality, the
increments $\delta v^{\nu,\epsilon,-,(k−1)}_i$, and $\delta v^{\nu,\epsilon,−,(k−1)}_{i,j}$
are extracted from (\ref{deltavk1}) straightforwardly.
Ad iii) we use the polynomial decay
\begin{equation}
\delta v^{\nu,\epsilon,-,(k)}_i(t,.)\in {\cal C}^4_{pol,m},~ k\geq 1,~ 1\leq i\leq D
\end{equation}
which holds also in the limit $k\uparrow \infty$, i.e., for all $1\leq i\leq D$ we have
\begin{equation}
\delta v^{\nu,\epsilon,-}_i(t,.)=\lim_{k\uparrow \infty }\delta v^{\nu,\epsilon,-,(k)}(t,.)\in {\cal C}^4_{pol,m}.
\end{equation}
This is enough polynomial decay at spatial infinity such that we can apply
compactness arguments in classical spaces (we do not even need Rellich’s compactness argument). We choose sequences $(\nu_l)_{l\geq 1}$ and $(\epsilon_l)_{l\geq 1}$ converging to zero
and consider the spatial transformation
\begin{equation}
v^{c,\nu_l,\epsilon_l,-}(t,y)=v^{\nu_l,\epsilon_l,-}(t,x)
\end{equation}
for $y_j=\arctan(x_j),~ 1\leq j\leq D$ and for all $t\in [0,T]$ for some time horizon $T>0$ such that local time contraction holds according to item ii), i.e., the argument of item ii) transfers. For multiindices $\gamma$ with $0\leq |\gamma|\leq m$, for all $t\in [0,T]$, for all $y\in (-2\pi,2\pi)^D$ and all $x\in {\mathbb R}^D$ we have
\begin{equation} 
{\big |}D^{\gamma}_y\delta_i^{c,\nu_l,\epsilon_l,-}(t,y){\big |}\leq c_0\left(1+|x|^{2m} \right){\big |}D^{\gamma}_y\delta_i^{\nu_l,\epsilon_l,-}(t,x){\big |}\leq C, 
\end{equation}
for some finite constants $c_0,C>0$. Then we apply
\begin{prop} For open and bounded $\Omega \subset {\mathbb R}^D$
and consider the function
space 
\begin{equation}
\begin{array}{ll}
C^m(\Omega):={\Big \{} f:\Omega\rightarrow {\mathbb R}|D^{\alpha}_xf~ \mbox{exists for $|\alpha|\leq m$}\\
\\
\hspace{2cm}\mbox{ and has a continuous extension to $\overline{\Omega}$}{\Big \}}
\end{array}
\end{equation}
where $\alpha = (\alpha_1, \cdots, \alpha_D)$ denotes a multiindex. Then the function space $C^m\left(\Omega \right)$ with the norm 
\begin{equation}
|f |_m := |f |_{C^m(\overline{\Omega})}=\sum_{|\alpha|\leq m}{\Big |}D^{\alpha}_xf{\Big |}
\end{equation}
is a Banach space. Here, 
\begin{equation}
{\big |}f{\big |}:=\sup_{x\in \Omega}{\big |}f(x){\big |}
\end{equation}
An analogous statement holds for $C^{m,\delta}(\Omega)$, where the derivatives up to order $m$ are H\"older continuous with exponent $\delta$.
\end{prop}
\section{Conclusion}
Solutions $v^E_i, 1\leq  i\leq  D = 3$ of the incompressible Euler Equation with regular
data $v^E_i(0,.)\in H^2\cap C^{1,\delta}$ which are singular at some time section $t = T$ on a thin
set of dimension 2 and in $H^2\cap C^{1,\delta}$
at all other times. They can exist next to
global classical solution branches for the same data. In the preceding example the first Euler velocity component $v^E_1$ is continuous while the first Euler vorticity
component $\omega^E_1$ is singular. The second and third velocity component are both singular at time $t = T$. In the example above we have $v^E_i(T,.)\in L^2\cap H^{1,1}$
for all $1\leq i \leq D$ at the critical time section at time $T$ , i.e., the first weak
derivatives of the velocity are in $L^1$. The first velocity component we have
$v^E_1(T,.)\in H^1\cap H^{2,1}$.
 This is a partial solution of a well-known problem described
in \cite{MB}. The constructed singular solutions of the Euler equation correspond to singular solutions $v^{\nu}_i,~ 1\leq i\leq D$ with viscosity $\nu >0$, and with time dependent force term
\begin{equation}\label{lap}
-\nu \Delta v^E_i
=: f_i, 1\leq i \leq D.
\end{equation}
Note that for the example above the force term $f_1(T,.)\in L^1$ and $f_i(T, .)\in H^{-2}\cap H^{-1,1}$ only. It is not possible to construct singular solutions of the Navier
Stokes equation with time dependent force terms in Schwartz space using the
construction above. On the other hand the construction above shows that Navier
stokes equations with time dependent force terms in $C^k\cap H^k$
for finite $k$ are not
necessarily smooth and note necessarily unique. The construction also shows
that uniqueness can be lost at any (short) time and that multiple solutions exist
for the incompressible Euler equation and corresponding solutions of the Navier
Stokes equation with related force terms defined in (\ref{lap}).

\end{document}